\documentclass{article}

\title{Two Symmetric Properties of Mersenne Numbers and Fermat Numbers}
\author{Shi Yongjin \thanks{\noindent\emph{Address}: Shenzhen, China. \emph{E-mail}: shi\_yongjin@163.com}}

\date{April, 2013}

\usepackage{amsmath}
\usepackage{amssymb}
\usepackage{amsthm}

\begin{document}

\maketitle
\begin{abstract}
Mersenne numbers and Fermat numbers are two hot and difficult issues in number theory. This paper constructs a special group for every positive odd number other than 1, and discovers an algorithm for determining the multiplicative order of $2$ modulo $q$ for each positive odd number $q$. It is worth mentioning that this paper discovers two symmetric properties of Mersenne numbers and Fermat numbers.\\

\noindent $\mathbf{Keywords}$: Mersenne numbers; Fermat numbers; Cycle; Commutative group; Algorithm; Multiplicative order; Symmetric properties
\end{abstract}

\section{Introduction}
A prime number is a positive integer $p>1$ that has no positive integer divisors other than 1 and  itself.
A Mersenne prime is a prime number of the form $M_n=2^n-1$. Primes of this form were first studied by Euclid who explored their relationship with the even perfect numbers. They are named after the French monk Marin Mersenne who studied them in the early 17th century.
It is easy to see that if $n$ is a composite number then so is $2^n-1$. In other words, in order for the Mersenne number $M_n$ to be prime, $n$ must be prime. The definition is therefore unchanged when written $M_p=2^p-1$ where $p$ is assumed prime.
More generally, numbers of the form $M_n=2^n-1$ without the primality requirement are called Mersenne numbers. Mersenne numbers are sometimes defined to have the additional requirement that $n$ be prime, equivalently that they be pernicious Mersenne numbers, namely those pernicious numbers whose binary representation contains no zeros. The smallest composite pernicious Mersenne number arises with $p = 11$.  This papers prefers to define a Mersenne number as a number of the form $M_p=2^p-1$ with $p$ restricted to prime values.
As of February 2013, 48 Mersenne primes are known. The largest known prime number $2^{57885161}-1$ is a Mersenne prime. Since 1997, all newly-found Mersenne primes have been discovered by the ¡°Great Internet Mersenne Prime Search¡± (GIMPS), a distributed computing project on the Internet. Many fundamental questions about Mersenne numbers and Mersenne primes remain unresolved. For example, are there infinitely many Mersenne primes? Is every Mersenne number $2^p-1$ square free? Can Zhou conjecture be established as true(see \cite{zhou})?\\

A Fermat number, named after Pierre de Fermat who first studied them, is a positive integer of the form $F_n=2^{2^n}+1$
where $n$ is a nonnegative integer. when a number of this form is prime, we call it a Fermat prime. The only known Fermat primes are the first five Fermat numbers: $F_0=3, F_1=5, F_2=17, F_3=257, F_4=65537$. Having discovered the fact that the first five Fermat numbers are all prime, Pierre de Fermat assumed that all numbers of this type were prime. But he was wrong. In 1732 after almost a century, Euler elegantly showed that $F_5=2^{2^5}+1=2^{32}+1=4294967297=641\times 6700417$. That year can be considered as the beginning of the search for divisors of other Fermat numbers. In 1770, Euler proved that every factor of $F_n$ must have the form $k\cdot2^{n+1}+1$. Edouard Lucas, improving the result by Euler, proved in 1878 that every factor of Fermat number, with $n$ at least 2, is of the form $k\cdot2^{n+2}+1$, where $k$ is a positive integer. This corollary is being used for discovery of Fermat number divisors. For 3 centuries more than 300 prime factors were found. So far only $F_0$ to $F_{11}$ have been completely factored. The distributed computing project Fermat Search is searching for new factors of Fermat numbers. Because of the scarcity and difficulty of finding these divisors, the person who discovers a new factor takes his place in history. Like Mersenne numbers and Mersenne primes, many fundamental questions about Fermat numbers and Fermat primes remain unresolved. For example, are there infinitely many Fermat primes? Is every Fermat number square free?\\

This paper tries to do some research into properties of Mersenne numbers and Fermat numbers.

\section{A special group for every positive odd number other than 1}

Let $q$ be any positive odd number larger than 1, we denote the set of all positive odd numbers less than $q$ by $G_{q}$, and denote the set of all positive odd numbers coprime to $q$ in $G_{q}$ by $G^\ast_q$, and define a binary operation on $G^\ast_q$ as follows:
$G^\ast_q\times G^\ast_q\longrightarrow G^\ast_q$ is given by $(a,b)\longmapsto \frac{ab-sq}{2^t}$ such that $0<ab-sq<q$ for some $s\in\mathbf{N}$ and $2\nmid\frac{ab-sq}{2^t}$ for some $t\in\mathbf{N}$. We denote $\frac{ab-sq}{2^t}$ by $a\ast b$, that is,  $a\ast b=\frac{ab-sq}{2^t}$.\\

Now we do an interesting calculation as follows:\\
1) Let $a_1$ be any odd number in $G^\ast_q$;\\
2) Let $a_2=\frac{q+a_1}{2^{k_1}}$ for some $k_1\in\mathbf{N}^*$ such that $2\nmid a_2$;\\
3) Let $a_3=\frac{q+a_2}{2^{k_2}}$ for some $k_2\in\mathbf{N}^*$ such that $2\nmid a_3$;\\
4) Repeat the above process, so that we can also get $a_4,a_5,\cdots$;\\
5) Untill some $n\in\mathbf{N}^*$ occurs such that $a_{n+1}=a_1$, we get a cycle $(a_1,a_2,\cdots,a_n)$, which can be called a cycle generated by $a_1$, absolutly it can also be generated by $a_2,\cdots,a_n$. $(a_1,a_2,\cdots,a_n)$ is also called a cycle of length $n$ or an $n$-cycle. Besides, we get a set $\{a_1,a_2,\cdots,a_n\}$, denote it by $\overline{a_1}$, we can also denote it by $\overline{a_i}$ with any $a_i\in\{a_1,a_2,\cdots,a_n\}$;\\
6) If $\overline{a_1}=G^\ast_q$, end the whole calculation. Otherwise we choose any positive odd number $b_1\in G^\ast_q-\overline{a_1}$, similarly we get a cycle $(b_1,b_2,\cdots,b_m)$. We also get a set $\{b_1,b_2,\cdots,b_m\}$, denote it by $\overline{b_k}(k=1,2,\cdots,m)$;\\
7) If $\overline{a_1}\cup \overline{b_1}=G^\ast_q$, end the whole calculation. Otherwise we choose any positive odd number $c_1\in G^\ast_q-\overline{a_1}-\overline{b_1}$, similarly we get a cycle $(c_1,c_2,\cdots,c_i)$. We also get a set $\{c_1,c_2,\cdots,c_i\}$, denote it by $\overline{c_r}(r=1,2,\cdots,i)$;\\
8) Repeat the process stated above untill $G^\ast_q=\overline{a_1}\cup\overline{b_1}\cup\overline{c_1}\cup\cdots\cup\overline{\alpha_1}\cup\overline{\beta_1}$ where $\overline{a_1},\overline{b_1},\overline{c_1},\cdots,\overline{\alpha_1},\overline{\beta_1}$ are all the sets we have got.\\
9) End the whole calculation. \\

If $\overline{x}=\overline{y}$, it is easy to show that there exists some $k\in\mathbf{N}$ such that $x\equiv 2^ky\pmod{q}$. Now we prove that if there exists some $k\in\mathbf{N}$ such that $x\equiv 2^ky\pmod{q}$, then $\overline{x}=\overline{y}$. If $x\equiv y\pmod{q}$, then $x=y$, thus $\overline{x}=\overline{y}$. If $x\equiv 2^{k_1}y\pmod{q}$ with $k_1>0$, there exists some positive odd number $r_1$ such that $x=2^{k_1}y-r_1q$. We denote $x$ by $\lambda_1$ and execute the calculation as above defined, then $q+\lambda_1=2^{k_1}y-(r_1-1)q$. If $r_1=1$, then $\lambda_2=y$, so $\overline{x}=\overline{y}$. If $r_1>1$, there exist some $s\in\mathbf{N^\ast}$ and some positive odd number $t$ such that $r_1-1=2^st$, then $q+\lambda_1=2^{k_1}y-2^stq$. Since $q+\lambda_1>0$, $y<tq$, then $k_1>s$, we have $\lambda_2=2^{k_1-s}y-tq$, denote $k_1-s$ by $k_2$ and denote $t$ by $r_2$, then $\lambda_2=2^{k_2}y-r_2q$. By analogy we have $\lambda_1=2^{k_1}y-r_1q,\lambda_2=2^{k_2}y-r_2q,\lambda_3=2^{k_3}y-r_3q,\cdots,\lambda_i=2^{k_i}y-r_iq,\cdots$. With the increase of $i$, positive odd number $r_i$ becomes smaller and smaller, hence there exists some $j\in\mathbf{N^\ast}$ such that $r_j=1$, then we have $\lambda_{j+1}=y$, therefore $\overline{x}=\overline{y}$.\\

We define a relation $R^\ast_q$ on $G^\ast_q\times G^\ast_q$ as follows:
$$R^\ast_q=\{(m,n)\in G^\ast_q\times G^\ast_q|\textrm{there exists }\overline{a_1}(a_1\in G^\ast_q)\textrm{ such that }\{m\}\cup\{n\}\subset\overline{a_1}\}.$$

It is easy to verify that the relation $R^\ast_q(\sim)$ on $G^\ast_q\times G^\ast_q$ is an equivalence relation on $G^\ast_q$ such that $x\sim y$ and $u\sim v$ imply $x\ast u\sim y\ast v$ for all $x,y,u,v\in G^\ast_q$. By convention $G^\ast_q/R^\ast_q$ denotes the set  of all equivalence classes of $G^\ast_q$ under $R^\ast_q$. We can define a binary operation on $G^\ast_q/R^\ast_q$ by $\overline{a}\overline{b}=\overline{a\ast b}$, where $\overline{x}$ denotes the equivalence class of $x\in G^\ast_q$.

First, we give the following proposition.

\newtheorem{proposition}{Proposition}
\begin{proposition}
For every odd number $q>1$, $G^\ast_q/R^\ast_q$ is a finite semigroup.
\end{proposition}
\begin{proof}
For all $x,y,z\in G^*_q$, $(\overline{x}\,\overline{y})\overline{z}=\overline{(x*y)*z}$, $\overline{x}(\overline{y}\,\overline{z})=\overline{x*(y*z)}$, there exist $s_1, s_2, s_3,s_4, t_1,t_2,t_3,t_4\in \mathbf{N}$ such that $x*y=\frac{xy-s_1q}{2^{t_1}}$,  $y*z=\frac{yz-s_3q}{2^{t_3}}$, $(x*y)*z=\frac{\frac{xy-s_1q}{2^{t_1}}z-s_2q}{2^{t_2}}$, $x*(y*z)=\frac{x\frac{yz-s_3q}{2^{t_3}}-s_4q}{2^{t_4}}$, where $0<xy-s_1q<q$, $0<yz-s_3q<q$, $0<\frac{xy-s_1q}{2^{t_1}}z-s_2q<q$, $0<x\frac{yz-s_3q}{2^{t_3}}-s_4q<q$, $2\nmid\frac{xy-s_1q}{2^{t_1}}$, $2\nmid\frac{yz-s_3q}{2^{t_3}}$, $2\nmid\frac{\frac{xy-s_1q}{2^{t_1}}z-s_2q}{2^{t_2}}$, $2\nmid\frac{x\frac{yz-s_3q}{2^{t_3}}-s_4q}{2^{t_4}}$.
Since $\frac{\frac{xy-s_1q}{2^{t_1}}z-s_2q}{2^{t_2}}\equiv 2^{\varphi(q)+t_3+t_4-t_1-t_2}\frac{x\frac{yz-s_3q}{2^{t_3}}-s_4q}{2^{t_4}}\pmod{q}$, where $\varphi(q)$ denotes Euler's totient  function that counts the number of positive integers not greater than and  relatively prime to $q$(see \cite{hardy}), it follows that $\overline{(x*y)*z}=\overline{x*(y*z)}$, then we have $(\overline{x}\,\overline{y})\overline{z}=\overline{x}(\overline{y}\,\overline{z})$ which means the binary operation defined by $\overline{a}\overline{b}=\overline{a*b}$ satisfies the associative property. Therefore,  $G^\ast_q/R^\ast_q$ is a finite semigroup.
\end{proof}

Next, we give the following theorem.

\newtheorem{theorem}{Theorem}
\begin{theorem}
For every odd number $q>1$, $G^\ast_q/R^\ast_q$ is a commutative group.
\end{theorem}

\begin{proof}
If $a,b,c,d\in G^\ast_q$, $\overline{a}=\overline{c}$ and $\overline{a}\overline{b}=\overline{c}\overline{d}$, then $\overline{a\ast b}=\overline{c\ast d}$. Since $a\ast b=\frac{ab-{s_1}q}{2^{t_1}}$ where $0<ab-{s_1}q<q$
for some $s_1\in\mathbf{N}$ and $\frac{ab-{s_1}q}{2^{t_1}}\equiv 1\pmod{2}$ for some $t_1\in\mathbf{N}$, and $c\ast d=\frac{cd-{s_2}q}{2^{t_2}}$ where $0<cd-{s_2}q<q$
for some $s_2\in\mathbf{N}$ and $\frac{cd-{s_2}q}{2^{t_2}}\equiv 1\pmod{2}$ for some $t_2\in\mathbf{N}$, there exists some $k\in\mathbf{N}$ such that $\frac{ab-{s_1}q}{2^{t_1}}\equiv 2^k \frac{cd-{s_2}q}{2^{t_2}}\pmod{q}$. We can choose some large enough number $k\in\mathbf{N}$ such that $k+t_1-t_2\geq 0$, thus $ab\equiv 2^{k+t_1-t_2}cd\pmod{q}$. Since $\overline{a}=\overline{c}$, there exists some $r\in\mathbf{N}$ such that $a\equiv 2^rc\pmod{q}$. We can also choose some large enough number $r\in\mathbf{N}$ such that $r+t_2-k-t_1\geq 0$, then $2^{r+t_2-k-t_1}b\equiv d\pmod{q}$, thus $\overline{b}=\overline{d}$. For all $a,b,c,d\in G^\ast_q$, $\overline{a}=\overline{c}$ and $\overline{a}\overline{b}=\overline{c}\overline{d}$ always imply $\overline{b}=\overline{d}$, hence the left cancellation law holds in $G^\ast_q/R^\ast_q$. Since $G^\ast_q/R^\ast_q$ is commutative, in which the right cancellation law also holds. The finite semigroup $G^\ast_q/R^\ast_q$ is commutative and cancellative, therefore $G^\ast_q/R^\ast_q$ is a  commutative group(see \cite{yang}).\\
\end{proof}

If $q$ is an odd composite number, for any $d\in G_q-G^{\ast}_q$, we can also obtain a cycle by the same technique as above defined. The cycle generated by $d$ can be called reducible cycle of $q$, conversely a cycle in $G^{\ast}_q/R^{\ast}_q$ is called irreducible cycle of $q$.

We can also define a relation $R_q$ on $G_q\times G_q$ as follows:
$$R_q=\{(m,n)\in G_q\times G_q|\textrm{there exists }\overline{a_1}(a_1\in G_q)\textrm{ such that }\{m\}\cup\{n\}\subset\overline{a_1}\}.$$

$G_q/R_q$ is defined to be the set of all equivalence classes of $G_q$ under $R_q$.
If we say all cycles of $q$, we mean all cycles in $G_q/R_q$; if we say all irreducible cycles of $q$, we mean all cycles in $G^\ast_q/R^\ast_q$; if we say all reducible cycles of $q$, we mean all cycles in $G_q/R_q-G^{\ast}_q/R^{\ast}_q$.\\

\section{An algorithm for determining the multiplicative order of $2$ modulo $q$ for each positive odd number $q$}

We denote the power of the prime factor 2 of an even number $n$ by $\tau(n)$. By using a method in \cite{hua}, we have\\
{\setlength\arraycolsep{2pt}
\begin{eqnarray*}
 \tau(\frac{(2q-2)!}{(q-1)!})&=&([\frac{2q-2}{2}]+ [\frac{2q-2}{2^2}]+[\frac{2q-2}{2^3}]+\cdots){}\\
 &&{}-([\frac{q-1}{2}]+ [\frac{q-1}{2^2}]+[\frac{q-1}{2^3}]+\cdots)\\
 &=&(q-1)+ ([\frac{q-1}{2}]+ [\frac{q-1}{2^2}]+[\frac{q-1}{2^3}]+\cdots){}\\
 &&{}-([\frac{q-1}{2}]+ [\frac{q-1}{2^2}]+[\frac{q-1}{2^3}]+\cdots)\\
 &=&q-1.
\end{eqnarray*}
Then we obtain $$\sum_{x\in G_q}\tau(q+x)=q-1.$$

For $x\in G^\ast_q$, we denote $\sum_{m\in \overline{x}}\tau(q+m)$ by $\xi(q,\overline{x})$.\\

Now we give an arithmetical property of group $G^\ast_q/R^\ast_q$ as follows.
\begin{proposition}
For every odd number $q>1$, let x be any odd number in $G^\ast_q$, then $\xi(q,\overline{x})$ is the smallest positive integer such that $q\mid 2^{\xi(q,\overline{x})}-1$.
\end{proposition}

\begin{proof}
Let $\alpha_1=x$, and execute the calculation as defined in Section 2, then we obtain a cycle $(\alpha_1,\alpha_2,\cdots,\alpha_n)$ generated by $x$. Let
$i_1,i_2,\cdots,i_n$ be positive integers such that $$q+\alpha_1=2^{i_1}\alpha_2,q+\alpha_2=2^{i_2}\alpha_3,\cdots,q+\alpha_n=2^{i_n}\alpha_1,$$
then
$$\alpha_1\equiv 2^{i_1}\alpha_2\pmod{q},\alpha_2\equiv 2^{i_2}\alpha_3\pmod{q},\cdots,\alpha_n\equiv 2^{i_n}\alpha_1\pmod{q},$$
thus $$\alpha_1\equiv 2^{i_1+i_2+\cdots+i_n}\alpha_1\pmod{q},$$
since $\alpha_1$ is coprime to $q$, then we have
$$2^{i_1+i_2+\cdots+i_n}\equiv 1\pmod{q}.$$

Now we show that for $v=1,2,3,\cdots,n$, $\alpha_1 \not\equiv 2^\sigma\alpha_v\pmod{q}$ with $0<\sigma<i_1$. It's clear $\alpha_1 \not\equiv 2^\sigma\alpha_2\pmod{q}$ with $0<\sigma<i_1$. Suppose there exist some $v\neq 2$ with  $1\leq v\leq n$ and some $\sigma$ with $0<\sigma<i_1$ such that $\alpha_1 \equiv 2^\sigma\alpha_v\pmod{q}$, then there exists some positive integer $l>1$ such that $lq+\alpha_1=2^\sigma\alpha_v$, since $q+\alpha_1=2^{i_1}\alpha_2$, we have $2^\sigma\alpha_v-2^{i_1}\alpha_2=(l-1)q$, then $\alpha_v-2^{i_1-\sigma}\alpha_2=\frac{l-1}{2^\sigma}q$, since $2^{i_1-\sigma}\alpha_2>0$, we have $\alpha_v>\frac{l-1}{2^\sigma}q\geq q$, it's a contradition, thus for $v=1,2,3,\cdots,n$, $\alpha_1 \not\equiv 2^\sigma\alpha_v\pmod{q}$ with $0<\sigma<i_1$. By analogy for $u,v=1,2,3,\cdots,n$, $\alpha_u \not\equiv 2^\sigma\alpha_v\pmod{q}$ with $0<\sigma<i_u$.

Suppose there exists some integer $r$ with $0<r<\sum_{k=1}^{n}i_k$ such that $2^{r}\equiv 1\pmod{q}$, then
\begin{equation}\label{eq:eps}
\alpha_1\equiv 2^{r}\alpha_1\pmod{q}.
\end{equation}
If $r=\sum_{k=1}^{j}i_k$ with $1\leq j< n$, since
\begin{equation}\label{eq:alp}
\alpha_1\equiv 2^{r}\alpha_{j+1}\pmod{q},
\end{equation}
from (\ref{eq:eps}) and (\ref{eq:alp}), we have
 $$\alpha_1 \equiv \alpha_{j+1}\pmod{q},$$
it's a contradiction, thus
$r\neq\sum_{k=1}^{j}i_k$ for $j=1,2,\cdots,n-1$. Then there exists some $h$ with $0\leq h\leq n-1$ such that $\sum_{k=0}^{h}i_k<r<\sum_{k=0}^{h+1}i_k(i_0=0)$, since $\alpha_1\equiv 2^{i_0+i_1+i_2+\cdots+i_h}\alpha_{h+1}\pmod{q}$, $\alpha_1\equiv 2^r\alpha_1\pmod{q}$, then $2^{i_0+i_1+i_2+\cdots+i_h}\alpha_{h+1}\equiv  2^r\alpha_1\pmod{q}$, thus $\alpha_{h+1}\equiv  2^{r-(i_0+i_1+i_2+\cdots+i_h)}\alpha_1\pmod{q}$ where $0<r-(i_0+i_1+i_2+\cdots+i_h)<i_{h+1}$, it's a contradiction, thus there exists no integer $h$ with $0\leq h\leq n-1$ such that $\sum_{k=0}^{h}i_k<r<\sum_{k=0}^{h+1}i_k$. It follows that there exists no integer $r$ with $0<r<\sum_{k=1}^{n}i_k$ such that $2^{r}\equiv 1\pmod{q}$.

Therefore, the proposition has been proved.

\end{proof}

Since $\xi(q,\overline{x})$ is the smallest positive integer such that $q\mid 2^{\xi(q,\overline{x})}-1$, then
for all $x,y\in G^\ast_q$,
$$\sum_{m\in \overline{x}}\tau(q+m)=\sum_{n\in \overline{y}}\tau(q+n),$$
that is,
\begin{equation}\label{eq:liu}
\xi(q,\overline{x})=\xi(q,\overline{y}).
\end{equation}
(\ref{eq:liu}) means the value of $\xi(q,\overline{x})$ is not dependent on $x$, but only dependent on $q$, hence we can denote $\xi(q,\overline{x})$ by $\varepsilon(q)$, which is a function of $q$ and denotes the multiplicative order of $2$ modulo $q$.\\

Proposition 2 can be expressed as the following form.
\begin{theorem}
For every odd number $q>1$, $\varepsilon(q)$ is the multiplicative order of $2$ modulo $q$.
\end{theorem}

Theorem 2 indeed gives a method of determining the multiplicative order of $2$ modulo $q$ for each positive odd number $q$.\\

Then we give a proposition as follows.
\begin{proposition}
Let $q$ be a positive odd number other than 1, then
$$|G^*_q/R^*_q|=\frac{\varphi(q)}{\varepsilon(q)}.$$
\end{proposition}

\newtheorem*{Corollary}{Corollary}
\begin{Corollary}

Let $p$ be a prime number other than 2, then

$$|G_p/R_p|=\frac{p-1}{\varepsilon(p)}.$$

\end{Corollary}

\section{Some symmetric properties of $G_q/R_q$ with $q$ of the form $2^n-1$ or $2^{2^n}+1$}

\subsection{A symmetric property of $G_q/R_q$ with $q$ of the form $2^{2^n}+1$}

For the first few Fermat numbers, we list all their cycles as follows.\\

\noindent $G_{3}/R_{3}:(1)$.\\

\noindent $G_{5}/R_{5}:(1,3)$.\\

\noindent $G_{17}/R_{17}:(1,9,13,15)$,  $(3,5,11,7)$.\\

\noindent $G_{257}/R_{257}$:\\
$(1,129,193,225,241,249,253,255)$, $(3,65,161,\,209,233,\,245,251,\,127)$,\\
$(5,\,131,\,97,\,177,\,217,\,237,\,247,63)$, $(7,33,145\,,201,229,\,243,125,\,191)$,\\
$(9,133,\,195,\,113,185,\,221,239,31)$, $(11,\,67,\,81,169,\,213,\,235,\,123\,,95)$,\\
$(13,135,\,49,153,\,205,231,\,61,159)$,$\;(15,17,137,197,227,121,189,223)$,\\
$(19,\,69,163,\,105,181,\,219,119,47)$,$\;(21,\:139,\:99,\:89,\,173,\:215,\:59,\:79)$,\\
$(23,\,35,\,73,\,165, 211, 117, 187, 111)$, $(25,\,141,199,57,157,\,207,\,29, 143)$,\\
$(27,\,71,\,41,\,149,\,203,\,115,\,93,\,\!175)$, $(37,\,147,101,\,179,109,183,55,\,39)$,\\
$(43,\;75,\;83,\;85,\;171,\:107,\;91,\;87)$, $\,(45,\,151,\,51,\,77,\,167,\,53,\,155,103)$.\\

\noindent $F_5=641\times6700417.$\\
Two cycles of 641 are listed as follows:\\
$(1,321,481,561,601,621,631,159,25,333,487,141,
391,129,385,513,577,609,\\625,\,633,637,639,\,5,323,241,441,541,\,591,77,359,125,383),$\\
$(3,161,401,521,581,611,313,477,559,75,179,205,423,133,387,257,449,545,\\
593,617,629,635,319,15,41,341,491,283,231,109,375,127).$\\
There are 32 numbers in each cycle.\\

From the above examples, this paper proposes a conjecture as follows.
\newtheorem{conjecture}{Conjecture}
\begin{conjecture}
For each natural number $n$, the cycles of $F_n$ are all $2^n$-cycles. Moreover, a positive integer $d$ other than 1 is factor of $F_n$ if and only if the cycles of $d$ are all $2^n$-cycles.
\end{conjecture}

\subsection{A symmetric property of $G_q/R_q$ with $q$ of the form $2^n-1$}

For the first few numbers of the form $2^n-1$, we list all their cycles as follows.\\
For $G_{2^n-1}/R_{2^n-1}$ and $G^\ast_{2^n-1}/R^\ast_{2^n-1}$, we'll observe there are how many 1-cycles in them, how many 2-cycles in them,
how many 3-cycles in them, $\cdots.$\\

\noindent $G_{3}/R_{3}:(1).$\\
$|G_{3}|=1.$\\

\noindent $G_{7}/R_{7}:(1);(3,5).$\\
$|G_{7}|=1\times 1+1\times 2=3.$\\

\noindent $G_{15}/R_{15}:(1),(5);(3,9);(7,11,13).$\\
$|G_{15}|=2\times1+1\times2+1\times3=7.$\\
$G^\ast_{15}/R^\ast_{15}:(1);(7,11,13).$\\
$|G^\ast_{15}|=1\times1+1\times3=4.$\\

\noindent $G_{31}/R_{31}:(1);(3,17),(5,9);(7,19,25),(11,21,13);(15,23,27,29).$\\
$|G_{31}|=1\times1+2\times2+2\times3+1\times4=15.$\\

\noindent $G_{63}/R_{63}:$\\
$(1),(9),(21);$\\
$(3,33),(5,17),(27,45);$\\
$(7,35,49),(11,37,25),(13,19,41);$\\
$(15,39,51,57),(23,43,53,29);$\\
$(31,47,55,59,61).$\\
$|G_{63}|=3\times1+3\times2+3\times3+2\times4+1\times5=31.$\\
$G^\ast_{63}/R^\ast_{63}:(1);(5,17);(11,37,25),(13,19,41);(23,43,53,29);(31,47,55,59,61).$\\
$|G^\ast_{63}|=1\times1+1\times2+2\times3+1\times4+1\times5=18.$\\

\noindent $G_{127}/R_{127}:$\\
$(1)$;\\
$(3,65),(5,33),(9,17);$\\
$(7,67,97),(11,69,49),(13,35,81),(19,73,25),(21,37,41);.$\\
$(15,71,99,113),(23,75,101,57),(27,77,51,89),(29,39,83,105),(43,85,53,45);$\\
$(31,79,103,115,121),(47,87,107,117,61),(55,91,109,59,93);$\\
$(63,95,111,119,123,125).$\\
$|G_{127}|=1\times1+3\times2+5\times3+5\times4+3\times5+1\times6=63.$\\

\noindent $|G_{2047}|=1\times1+5\times2+15\times3+30\times4+42\times5+42\times6+30\times7+15\times8+5\times9+1\times10=1023.$\\
$|G^\ast_{2047}|=1\times1+5\times2+14\times3+28\times4+40\times5+40\times6+28\times7+14\times8+5\times9+1\times10=968.$\\
$M_{11}=2047=23\times89.$\\
$G_{23}/R_{23}:(1,3,13,9);(5,7,15,19,21,11,17).\\$
$|G_{23}|=1\times4+1\times7=11.$\\
$G_{89}/R_{89}$:\\
$(3,23,7);\\
(1,45,67,39);\\
(5,47,17,53,71),(13,51,35,31,15);\\
(9,49,69,79,21,55),(19,27,29,59,37,63);\\
(11,25,57,73,81,85,87);\\
(33,61,75,41,65,77,83,43).\\$
$|G_{89}|=1\times3+1\times4+2\times5+2\times6+1\times7+1\times8=44.$\\

From the above examples, we obtain the following expansions:\\
$\frac{M_3-1}{2}=2^2-1=3=1\times 1+1\times 2.$\\
$\frac{M_5-1}{2}=2^4-1=15=1\times1+2\times2+2\times3+1\times4.$\\
$\frac{M_7-1}{2}=2^6-1=63=1\times1+3\times2+5\times3+5\times4+3\times5+1\times6.$\\
$\frac{M_{11}-1}{2}=2^{10}-1=1023=1\times1+5\times2+15\times3+30\times4+42\times5+42\times6+30\times7+15\times8+5\times9+1\times10.$\\
$\frac{M_{13}-1}{2}=2^{12}-1=4095=1\times1+6\times2+22\times3+55\times4+99\times5+132\times6+132\times7+99\times8+55\times9+22\times10+6\times11+1\times12.$\\
$\cdots$\\

For each odd prime number $p$,\\
$M_p-1=2^p-2=(1+1)^p-2=2\binom{p}{1}+2\binom{p}{2}+2\binom{p}{3}+\cdots+2\binom{p}{\frac{p-1}{2}},$\\
{\setlength\arraycolsep{1pt}
\begin{eqnarray}\label{eq:jiao}
\frac{M_p-1}{2}&=&\binom{p}{1}+\binom{p}{2}+\binom{p}{3}+\cdots+\binom{p}{\frac{p-1}{2}}\nonumber\\
&=&\frac{(p-1)!}{1!(p-1)!}(1+p-1)+\frac{(p-1)!}{2!(p-2)!}(2+p-2)+\frac{(p-1)!}{3!(p-3)!}(3+p-3){}\nonumber\\
&&{}+\cdots+\frac{(p-1)!}{(\frac{p-1}{2})!(\frac{p+1}{2})!}(\frac{p-1}{2}+\frac{p+1}{2})\nonumber\\
&=&\frac{(p-1)!}{1!(p-1)!}\cdot 1+\frac{(p-1)!}{2!(p-2)!}\cdot 2+\cdots+\frac{(p-1)!}{(\frac{p-1}{2})!(\frac{p+1}{2})!}\cdot\frac{p-1}{2}{}\nonumber\\
&&{}+\frac{(p-1)!}{(\frac{p+1}{2})!(\frac{p-1}{2})!}\cdot\frac{p+1}{2}+\cdots+\frac{(p-1)!}{(p-2)!2!}\cdot(p-2)+\frac{(p-1)!}{(p-1)!1!}\cdot(p-1)\nonumber\\
&=&\sum_{k=1}^{p-1}\frac{(p-1)!}{k!(p-k)!}k.\nonumber\\
\end{eqnarray}

From (\ref{eq:jiao}) and the number of cycles of different length in $G_3/R_3$,$G_7/R_7$,$G_{31}/R_{31}$,\\
$G_{127}/R_{127}$,$G_{2047}/R_{2047}$,$G_{23}/R_{23}$ and $G_{89}/R_{89}$, this paper proposes a conjecture as follows.
\begin{conjecture}
For each prime number $p$, the number of $k$-cycles of $M_p$ is $\frac{(p-1)!}{k!(p-k)!}$. Moreover, a positive integer $d$ other than 1 is factor of $M_p$ if and only if the number of $k$-cycles of $d$ is equal to the number of $(p-k)$-cycles of $d$.
\end{conjecture}

\section{Conclusions}
This paper obtains four main conclusions as follows:\\
1) For every odd number $q>1$, $G^\ast_q/R^\ast_q$ is a commutative group.\\
2) For every odd number $q>1$, $\varepsilon(q)$ is the multiplicative order of $2$ modulo $q$.\\
3) For each natural number $n$, the cycles of $F_n$ are all $2^n$-cycles. Moreover, a positive integer $d$ other than 1 is factor of $F_n$ if and only if the cycles of $d$ are all $2^n$-cycles.\\
4) For each prime number $p$, the number of $k$-cycles of $M_p$ is $\frac{(p-1)!}{k!(p-k)!}$. Moreover, a positive integer $d$ other than 1 is factor of $M_p$ if and only if the number of $k$-cycles of $d$ is equal to the number of $(p-k)$-cycles of $d$.\\

However, this paper doesn't give the proofs of 3) and 4), so the reader who is interested in the two problems could do further research.


\begin{thebibliography}{99}
\bibitem{zhou}Zhou Haizhong. \emph{The distribution of Mersenne primes}. Acta Scientiarum Naturalium Universitatis Sunyatseni, 1992, 31(4): 121-122.
\bibitem{hardy}G. H. Hardy, E. M. Wright. \emph{An introduction to the theory of numbers}. Beijing: Posts \& Telecom Press, 2009: 63-65.

\bibitem{yang}Yang Zixu. \emph{Modern algebra, second edition}. Beijing: High Education Press, 2003: 37-38.

\bibitem{hua}Hua Luogeng. \emph{An introduction to the theory of numbers}. Beijing: Science Press, 1957: 15-16.
\end{thebibliography}
\end{document}